\documentclass{article}

\newtheorem{theorem}{Theorem}

\newtheorem{lemma}[theorem]{Lemma}

\usepackage{latexsym}

\title{A constructive proof of a theorem by Ferreira-Zantema
}
\author{Toshiyasu Arai
\\
Graduate School of Mathematical Sciences
\\
University of Tokyo
\\
3-8-1 Komaba, Meguro-ku,
Tokyo 153-8914, JAPAN
\\
tosarai@ms.u-tokyo.ac.jp
}
\date{}
\begin{document}
\maketitle

\begin{abstract}
This note was written in Jan.\,23, 2015 to answer a problem raised by G. Moser,
who asked a constructive proof of a theorem by Ferreira-Zantema\cite{FerreiraZantema}.
\end{abstract}

For a binary relation $<$ on a set $T$, $W(<)$ denotes the well-founded part of $T$ with respect to $<$.

Let $\mathcal{F}$ be a non-empty and finite set of function symbols with $\#\mathcal{F}=N>0$.
Each $f$ has a fixed arity $ar(f)\in\omega$.
When writing $f(a)$, we tacitly assume that $a\in T(\mathcal{F})^{ar(f)}$.
If $a=(a_{0},\ldots,a_{n-1})$, then
$t\in a$ iff $t=a_{i}$ for an $i<n$.

Let $<$ be a proper order, i.e., irreflexive and transitive relation on $T(\mathcal{F})$, and for each $f\in\mathcal{F}$,
a relation (not necessarily a proper order) $<_{f}$ on $T(\mathcal{F})^{ar(f)}$ is given.

Assume the following three conditions for $f\in\mathcal{F}$:
\begin{enumerate}
\item
$<$ contains the subterm relation: for any proper subterm $s$ of $t$, $s<t$.

\item
If $f(b)<f(a)$, then either $f(b)\leq a_{i}$ for some $a_{i}\in a$, or
$b<_{f}a$.

\item
\begin{equation}\label{eq:<f}
\forall a\in T(\mathcal{F})^{ar(f)}[a\subset W(<) \Rightarrow a\in W(<_{f})]
\end{equation}

\end{enumerate}

\begin{theorem}\label{thm2.2}{\rm (Ferreira-Zantema\cite{FerreiraZantema})}
\\
$T(\mathcal{F})\subset W(<)$.
\end{theorem}

We show the following lemma first.

\begin{lemma}\label{lemma:a}
Let $\{f_{1},\ldots,f_{K}\}\subset\mathcal{F}$ be a set of distinct function symbols, and $a_{i}\in T(\mathcal{F})^{ar(f_{i})}$.
Let for $i=0,1,\ldots,K$,
\begin{eqnarray*}
T_{f,0} & := & T(\mathcal{F})^{ar(f)}
\\
T_{f,i} & := & \{b\in T(\mathcal{F})^{ar(f)}: f(b)<f_{i}(a_{i})\}\, (i\neq 0)
\\
A_{f,i}(b) & :\Leftrightarrow & [b\in T_{f,i} \to b\subset W(<) \to f(b)\in W(<)]
\\
IH_{f,i}(a) & :\Leftrightarrow & \forall b<_{f}a\, A_{f,i}(b)
\end{eqnarray*}

Assume the following three for any $i=1,\ldots,K$:
\begin{enumerate}
\item\label{lemma:a1}
$f_{i+1}(a_{i+1})<f_{i}(a_{i})$.
\item\label{lemma:a2}
$IH_{f_{i},i-1}(a_{i})$.
\item\label{lemma:a3}
$a_{i}\subset W(<)$.
\end{enumerate}
Then $f_{K}(a_{K})\in W(<)$.
\end{lemma}
{\bf Proof} by (meta)induction on $N-K$.

By MIH we have for any $f\not\in\{f_{1},\ldots,f_{K}\}$,
$\forall b[f(b)<f_{K}(a_{K}) \to IH_{f,K}(b)\to b\subset W(<) \to f(b)\in W(<)]$, i.e.,
$\forall b[\forall c<_{f}b\, A_{f,K}(c) \to A_{f,K}(b)]$.
In other words, the $\Pi^{1}_{2}$-predicate $A_{f,K}$ is progressive with respect to $<_{f}$.
Hence $\Pi^{1}_{2}\mbox{-BI}_{0}\vdash W(<_{f})\subset A_{f,K}$.

We show that
\[
s<f_{K}(a_{K}) \to s\in W(<)
\]
by subsidiary induction on the size of terms $s$.

Let $s=f(b)<f_{K}(a_{K})$.
For any $b_{i}\in b$, $b_{i}<f_{K}(a_{K})$ since $<$ contains the subterm relation and transitive.
By SIH we have $b\subset W(<)$.

First consider the case $f=f_{i}$ for some $1\leq i\leq K$.
Then  by (\ref{lemma:a1}) and the transitivity of $<$, we have
$s=f_{i}(b)<f_{i}(a_{i})$.

If $s\leq a_{ij}$ for an $a_{ij}\in a_{i}$, then (\ref{lemma:a3}) yields $s\leq a_{ij}\in W(<)$, and $s\in W(<)$.
Otherwise we have $b<_{f}a_{i}$.
$IH_{f_{i}}(a_{i})$, (\ref{lemma:a2}), with  $b\subset W(<)$ yields $s=f_{i}(b)\in W(<)$.

Second consider the case $f\not\in\{f_{1},\ldots,f_{K}\}$.
By MIH we have $W(<_{f})\subset A_{f,K}$.
$b\subset W(<)$ with (\ref{eq:<f}) yields $b\in W(<_{f})$, and $A_{f,K}(b)$.
On the other hand we have $f(b)<f_{K}(a_{K})$, i.e., $b\in T_{f,K}$.
Thus we conclude $f(b)\in W(<)$.

\hspace*{\fill} $\Box$
\\

From Lemma \ref{lemma:a} we see that for any $f\in\mathcal{F}$,
$\forall a[IH_{f,0}(a) \to a\subset W(<) \to f(a)\in W(<)]$.
In other words, the $\Pi^{1}_{2}$-predicate $A_{f,0}$ is progressive with respect to $<_{f}$.
Hence $(\Pi^{1}_{2}\mbox{-BI}_{0}\vdash) W(<_{f})\subset A_{f,0}$, i.e.,
$a\in W(<_{f}) \to a\subset W(<) \to f(a)\in W(<)$.
Under (\ref{eq:<f}) this is equivalent to
\[
\forall a[ a\subset W(<) \to f(a)\in W(<)]
\]
Now Theorem \ref{thm2.2}, $t\in W(<)$ is seen by induction on the size of terms $t$.
Thus we get the following Theorem \ref{th:each}.

\begin{theorem}\label{th:each}
Theorem \ref{thm2.2} is proved in $\Pi^{1}_{2}\mbox{{\rm -BI}}_{0}$ for
{\rm each} finite $\mathcal{F}$.
\end{theorem}

When $\mathcal{F}$ ranges over finite sets of function symbols in Theorem \ref{thm2.2},
we need to prove Lemma \ref{lemma:a} by formal complete induction, and
the lemma is a sentence of the form $\forall n\exists X\forall Y\, \theta$ with a first-order $\theta$.
Therefore the proof of the lemma is formalizable in $\Pi^{1}_{2}\mbox{-BI}$ with full induction schema, but not in 
$\Pi^{1}_{2}\mbox{-BI}_{0}$ with restricted induction.

Indeed, Theorem \ref{thm2.2} for \textit{any} finite set $\mathcal{F}$ implies the well-foundedness 
$W(\vartheta(\Omega^{\omega}))$ of the proof-theoretic
ordinal $\vartheta(\Omega^{\omega})$ of $\Pi^{1}_{2}\mbox{-BI}_{0}$(Lemma \ref{lem:wo} below), and is not provable in $\Pi^{1}_{2}\mbox{-BI}_{0}$.
It is known that Kruskul's theorem is equivalent to $W(\vartheta(\Omega^{\omega}))$ over $\mbox{ACA}_{0}$, cf.  \cite{RathjenWeiermann}.
Therefore we obtain the following Theorem \ref{th:equivalenceany}.

\begin{theorem}\label{th:equivalenceany}
Over $\mbox{{\rm ACA}}_{0}$ the following facts are equivalent each other:
\begin{enumerate}
\item
Any simplification order over {\rm any} finite $\mathcal{F}$ is terminating.
\item
Theorem \ref{thm2.2} for {\rm any} finite set $\mathcal{F}$.
\item
$W(\vartheta(\Omega^{\omega}))$.
\item
Kruskul's theorem $KT(\omega)$ for {\rm any} finite trees.
\end{enumerate}
\end{theorem}

\begin{theorem}\label{th:equivalenceeach}
Over $\mbox{{\rm ACA}}_{0}$ the following facts are equivalent each other:
\begin{enumerate}
\item
Any simplification order over {\rm each} finite $\mathcal{F}$ is terminating.
\item
Theorem \ref{thm2.2} for {\rm each} finite set $\mathcal{F}$.
\item
$W(\vartheta(\Omega^{k}))$ for $k=0,1,2,\ldots$
\item
Kruskul's theorems $KT(k)$ for $k$-branching trees for $k=0,1,2,\ldots$
\end{enumerate}
\end{theorem}

\begin{lemma}\label{lem:wo}
Theorem \ref{thm2.2} for \textit{any} finite set $\mathcal{F}$ implies the well-foundedness 
$W(\vartheta(\Omega^{\omega}))$.
\end{lemma}
{\bf Proof}.

Each ordinal$<\vartheta(\Omega^{k+1})$ is represented by a term over the symbols $0,1,\ldots,k,+,\vartheta$ and $\Omega$.
Let $\mathcal{F}_{k}=\{f_{i}: i\leq k\}\cup\{g,1\}$,
where $ar(f_{i})=i+1$, $ar(g)=2$ and $ar(1)=0$.
Each ground term $t\in T(\mathcal{F})$ denotes a non-zero ordinal $o(t)<\vartheta(\Omega^{k+1})$ as follows.
$o(1)=1$, $o(g(t,s))=o(t)\# o(s)$ for the natural sum $\alpha\#\beta$ of ordinals $\alpha$ and $\beta$, and
$o(f_{i}(t_{i},\ldots,t_{0}))=\vartheta(\Omega^{i}o(t_{i})+\cdots+\Omega^{0}o(t_{0}))$.
Let
\[
t\prec s :\Leftrightarrow o(t)<o(s) 
.\]
The relation $\prec$ on terms is a proper order, and contains the subterm relation.
$\alpha=\vartheta(\Omega^{i}\alpha_{i}+\cdots+\Omega^{0}\alpha_{0})<
\vartheta(\Omega^{i}\beta_{i}+\cdots+\Omega^{0}\beta_{0})=\beta$ iff
either $\alpha\leq\beta_{j}$ for a $j$, or $\alpha_{j}<\beta$ for any $j$ and 
$(\alpha_{i},\ldots,\alpha_{0})<_{lx}(\beta_{i},\ldots,\beta_{0})$ for the lexicographic order $<_{lx}$.
Therefore if $f_{i}(t_{i},\ldots,t_{0})\prec f_{i}(s_{i},\ldots,s_{0})$, then either $f_{i}(t_{i},\ldots,t_{0})\preceq s_{j}$ for a $j$, or
$(t_{i},\ldots,t_{0})\prec_{lx}(s_{i},\ldots,s_{0})$, where $\prec_{lx}$ is the lexicographic order on tuples of the same lengths
induced by $\prec$.
Moreover if $g(t_{0},t_{1})\prec g(s_{0},s_{1})$, then 
$(t_{0},t_{1})\prec_{m}(s_{0},s_{1}):\Leftrightarrow \exists i,j<2(t_{i}\prec s_{j} \land t_{1-i}\preceq s_{j-1})$, 
where $\prec_{m}$ is a multiset extension of $\prec$.
Hence (\ref{eq:<f}) is enjoyed for each function symbol.
Theorem \ref{thm2.2} yields $T(\mathcal{F}_{k})\subset W(\prec)$.

For each ordinal (term in normal form) $\alpha<\vartheta(\Omega^{k+1})$, let $\alpha^{+}$ denote the ordinal
obtained from $\alpha$ by replacing the subterms $0$ by $1$.
$0^{+}=1$, $(\alpha\#\beta)^{+}=\alpha^{+}\#\beta^{+}$ and 
$(\vartheta(\Omega^{i}\alpha_{i}+\cdots+\Omega^{0}\alpha_{0}))^{+}=
\vartheta(\Omega^{i}\alpha_{i}^{+}+\cdots+\Omega^{0}\alpha_{0}^{+})$.
Then it is easy to see that $\alpha<\beta$ iff $\alpha^{+}<\beta^{+}$.
Moreover for each $\alpha<\vartheta(\Omega^{k})$ there exists a term $t$ such that $o(t)=\alpha^{+}$.
Therefore $T(\mathcal{F}_{k})\subset W(\prec)$ yields $W(\vartheta(\Omega^{k+1}))$.

\hspace*{\fill} $\Box$

\end{document}